\documentclass[12pt]{amsart}

\usepackage{amsthm}
\usepackage{enumerate}

\usepackage[colorlinks]{hyperref}

\newtheorem{thm}{Theorem}[section]
\newtheorem{lem}[thm]{Lemma}
\newtheorem{prop}[thm]{Proposition}
\newtheorem{cor}[thm]{Corollary}
\newtheorem{rem}[thm]{Remark}

\title[Spectral Synthesis in the Multidimensional Fourier Algebra]{Spectral Synthesis in the Multidimensional Fourier Algebra and the Varopolous Algebra for Compact Groups}
\author{Kanupriya}
\address{Kanupriya,\newline\indent Department of Mathematics,\newline\indent Indian Institute of Technology Delhi,\newline\indent New Delhi - 110016, India.}
\email{kanupriyawadhawan3@gmail.com}
\author{N. Shravan Kumar}
\address{N. Shravan Kumar,\newline\indent Department of Mathematics,\newline\indent Indian Institute of Technology Delhi,\newline\indent New Delhi - 110016, India.}
\email{shravankumar.nageswaran@gmail.com}

\begin{document}

\begin{abstract}
	Let $G$ be a compact group and let $A^n(G)$ denote the multidimensional Fourier algebra introduced by Todorov and Turowska. In this note, we first define the multidimensional version of the Varopolous algebra and show that the multidimensional Fourier algebra can be embedded into the multidimensional Varopolous algebra. Using this embedding, we also prove a result on parallel synthesis, subsuming the earlier results of Varopolous, Spronk-Turowska and Parthsarathy-Prakash.
\end{abstract}
	
	

\maketitle

\section{Introduction}
	Spectral synthesis in locally compact abelian groups is a vintage topic in abelian harmonic analysis. Among all the classical theorems, the most celebrated one is the Malliavin's theorem on the failure of spectral synthesis in non-discrete groups, going back to 1959. Varopolous \cite{Var} pioneered tensor product methods to establish Malliavin's theorem. His proof is based on a relation that he established between spectral synthesis in the Fourier algebra and the projective tensor product of a commutative $C^\ast$-algebra with itself. More specifically, if $G$ is a compact abelian group, there is a natural embedding of $A(G)$ into $V(G)=C(G)\hat{\otimes}C(G).$ Using this embedding, he showed that, if $E$ is a closed subset of $G,$ then $E$ is a set of synthesis for $A(G),$ if and only if $E^\ast=\{(x,y)\in G\times G:x+y\in G\}$ is a set of synthesis for $V(G).$ This tensor product algebra $V(G)$ was later called as the Varopolous algebra. 
		
	Let $G$ be a compact group. Then, on $C(G)\otimes C(G),$ the Haagerup norm is exactly Grothendieck’s H-norm, which in turn, by Grothendieck’s inequality, is equivalent to the projective tensor norm. Thus, in 2002, Spronk and Turowska replaced the projective tensor product by Haagerup tensor product. With this as the definition of the Varopolous algebra, they were able to extend the result of Varopolous to the context of compact groups. For an extension of the result of Spronk and Turowska, see \cite{PS1}.
		
	The multidimensional Fourier algebra $A^n(G)$ on a locally compact group $G$ was introduced a decade ago. For the case when $n=2$ and also when $G$ is abelian, this space, via the Fourier transform, concides with the space of bimeasures, introduced by Graham and Schreiber \cite{GS}. This paper studies spectral synthesis for this algebra. More specifically, we discuss the multidimensional analogue of the parallel synthesis result on compact groups. 
		
	In section \ref{MVA}, we introduce the multidimensional analogue of the Varopolous algebra and show that it is a commutative Banach algebra. In section \ref{PNMap}, we define the analogues of the maps $P$ and $N.$ Using this we provide an embedding of the multidimensional Fourier algebra into the multidimensional Varopolous algebra. See \cite{Blei} for another extension of these maps, in the case of compact abelian groups. Finally, in section \ref{SSAV}, we prove the relation between spectral synthesis in multidimensional Fourier algebra and the multidimensional Varopolous algebra. 
		
	In \cite{KL1}, the authors have defined the notion of $X$-synthesis, extending the classical notion of synthesis. This was taken up and studied further by Parthasarathy and Prakash \cite{PP1,PP2}. In \cite{PP1}, the authors show how the result on parallel synthesis can be extended to $X$-synthesis. Theorem \ref{SSFAVA} of section \ref{SSAV}, gives the result of parallel synthesis for the multidimensional case. We would like to mention here that this result is new even when the underlying group is abelian. Some preparations for this result are worked out in section \ref{SMAV}.
		
	It is a natural question to ask the following question: {\it does the above mentioned result on parallel synthesis extend to Ditkin sets?} This question was taken up in \cite{PS2} and was answered affirmatively in the case of strong Ditkin sets. In Theorem \ref{DFAVA}, we take up this question for the multidimensional case, but for the strong Ditkin sets.
		
	We shall now begin with some preliminaries that are needed in the sequel.

\section{Preliminaries}

    \subsection{Fourier analysis on compact groups}
        Let $G$ be a compact group. Then it is well known that $G$ possesses a unique Haar measure $dx$ such that $\underset{G}{\int}dx=1.$ Further, an irreducible unitary representation of $G$ is always finite-dimensional. Let $\widehat{G}$ denote the set of all irreducible unitary representations of $G.$ Then $\widehat{G}$ is called the unitary dual of $G$ We shall denote by $d_\pi$ the dimension of the representation $\pi.$ 
        
        Given a representation $\pi$ and $u,v\in \mathcal{H}_\pi,$ the mapping $x\mapsto<\pi(x)u,v>$ is called the coefficient function of $\pi.$ If $\pi\in\widehat{G},$ let $\mathcal{E}_\pi$ denote the subspace of $L^2(G)$ spanned by the matrix coefficients of $\pi.$ Note that the space $\mathcal{E}_\pi$ is of dimension $d_\pi^2$ and so let $\{e_i^\pi:i=1,2,\ldots,d_\pi^2\}$ denote an orthonormal basis for $\mathcal{E}_\pi.$ Next observe that each $\mathcal{E}_\pi$ is invariant under $\lambda.$ 
        \begin{thm}
        Let $G$ be a compact group.
        \begin{enumerate}[(i)]
        \item The coefficient function arising out of an irreducible unitary representation belongs to $L^2(G).$
        \item{\rm(Peter-Weyl Theorem)} The space $L^2(G)$ is equal to the Hilbert space direct sum of coefficient spaces of the irreducible unitary representations of $G,$ i.e., $$L^2(G)=\underset{\pi\in\widehat{G}}{\bigoplus}\ \mathcal{E}_\pi.$$
        \end{enumerate}
        \end{thm}
        For any finite subset $F$ of $\widehat{G},$ let $P_F$ denote the orthogonal projection from $L^2(G)$ onto $\mathcal{E}_F:=\underset{\pi\in F}{\bigoplus}\ \mathcal{E}_\pi.$ We shall denote by $\lambda,$ the left regular representation of $G$ on $L^2(G),$ given by $$\lambda(x)(f)(y):=f(x^{-1}y).$$ Note that $P_F$ commutes with the left regular representation, i.e., $P_F\lambda(x)=\lambda(x)P_F$ for all $x\in G.$ Also, every irreducible unitary representation of a compact group can be realized as a subrepresentation of $\lambda.$
        
	\subsection{Operator spaces and their tensor products}
        In this section we review some basics on operator spaces. As tensor products will play a major role, our major aim of this section is to collect the required preliminaries on these topics.

        Let $X$ be a linear space. By $M_n(X)$ we shall mean the space of all $n\times n$ matrices with entries from the space $X.$ An {\it operator space} is a complex vector space $X$ together with an assignment of a norm $\|\cdot\|_n$ on the matrix space $M_n(X),$ for each $n\in\mathbb{N},$ such that
        \begin{enumerate}[(i)]
            \item $\|x\oplus y\|_{m+n}=max\{\|x\|_m,\|y\|_n\}$ and
            \item $\|\alpha x\beta\|_n\leq\|\alpha\|\|x\|_m\|\beta\|$
        \end{enumerate}
        for all $x\in M_m(X),$ $y\in M_n(X),$ $\alpha\in M_{n,m}$ and $\beta\in M_{m,n}.$ 

        Let $X$ and $Y$ be operator spaces and let $\varphi:X\rightarrow Y$ be a linear transformation. For any $n\in\mathbb{N},$ the {\it $n^{\mbox{th}}$-amplification} of $\varphi,$ denoted $\varphi_n,$ is defined as a linear transformation $\varphi_n:M_n(X) \rightarrow M_n(Y)$ given by $\varphi_n([x_{ij}]):=[\varphi(x_{ij})].$ The linear transformation $\varphi$  is said to be {\it completely bounded} if $\sup\{\|\varphi_n\|:n\in\mathbb{N}\}<\infty.$ We shall denote by $\mathcal{CB}(X,Y)$ the space of all completely bounded linear mappings from $X$ to $Y$ equipped with the norm, denoted $\|\cdot\|_{cb},$ $$\|\varphi\|_{cb}:=\sup\{\|\varphi_n\|:n\in\mathbb{N}\},\ \varphi\in\mathcal{CB}(X,Y).$$ We shall say that $\varphi$ is a {\it complete isometry} if $\varphi_n$ is an isometry $\forall\ n\in\mathbb{N}.$

        By Ruan's theorem, given an abstract operator space $X$ there exists a Hilbert space $\mathcal{H}$ and a closed subspace $Y\subseteq\mathcal{B}(\mathcal{H})$ such that $X$ and $Y$ are completely isometric. 

        Given operator spaces $X$ and $Y$ the {\it Haagerup tensor norm} of $x\in M_n(X\otimes Y)$ is given as $$\|u\|_h=\inf\left\{ \|x\|\|y\|:u=x\odot y,x\in M_{n,r}(X), y\in M_{r,n}(Y),r\in\mathbb{N} \right\}.$$ The quantity $\|\cdot\|_h$ is an operator space norm and we shall denote by $X\otimes^h Y$ the resulting operator space. In case $X$ and $Y$ are C*-algebras, for the case $n=1,$ the Haagerup norm can be written as follows. For $u\in X\otimes^h Y,$ $$\|u\|_h=\inf \left\{ \left\| \underset{n\in\mathbb{N}}{\sum} x_n x_n^\ast \right\|_X^{1/2} \left\| \underset{n\in\mathbb{N}}{\sum} y_n^\ast y_n \right\|_Y^{1/2} :u= \underset{n\in\mathbb{N}} {\sum}x_n\otimes y_n \right\}.$$ The following is the multidimensional analogue of \cite[Theorem 4.3]{Sm}. As the proof of this follows the same lines as in \cite{Sm} we skip it.
        \begin{thm}\label{WASm}
            Let $\mathcal{H}$ be a Hilbert space and let $v=v_1\otimes v_2 \otimes \ldots \otimes v_{n+1} \in \otimes_{n+1}^h\mathcal{B}(\mathcal{H}).$ Define $\phi_v:\otimes_{n}^h\mathcal{K}(\mathcal{H})\rightarrow\mathcal{K}(\mathcal{H})$ as $\phi(k)=v_1 k_1 v_2 k_2 v_3\ldots v_n k_n v_{n+1},$ where $k=k_1\otimes k_2\otimes\ldots\otimes k_n.$ Then the mapping $v\mapsto\phi_v$ from $\otimes_{n+1}^h\mathcal{B}(\mathcal{H})$ to $\mathcal{CB}(\otimes_{n}^h\mathcal{K}(\mathcal{H}), \mathcal{K}(\mathcal{H}))$ is an isometry.
        \end{thm}
            
        The {\it extended Haagerup tensor product} of $X$ and $Y,$ denoted $X\otimes^{eh} Y$ is defined as the space of all normal multiplicatively bounded functionals on $X^\ast\times Y^\ast.$ By \cite{ER2}, $(X\otimes^h Y)^\ast$ and $X^\ast\otimes^{eh}Y^\ast$ are completely isometric. Given dual operator spaces $X^\ast$ and $Y^\ast,$ the {\it $\sigma$-Haagerup tensor product (or normal Haagerup tensor product)} is defined by $$X^\ast\otimes^{\sigma h}Y^\ast=(X\otimes^{eh}Y)^\ast.$$ Further, the following inclusions hold completely isometrically: $$X^\ast\otimes^h Y^\ast \hookrightarrow X^\ast\otimes^{eh} Y^\ast \hookrightarrow X^\ast\otimes^{\sigma h} Y^\ast.$$

        For more on these tensor products, we refer to \cite{ER1,ER2}. For further reading on operator spaces, the reader is asked to refer \cite{ER3} or \cite{Pis}.

    \subsection{Multidimensional Fourier algebra}
		We first give a short introduction to the Fourier algebra.
		
		Let $G$ be a locally compact group. Let $\pi$ be a unitary representation of $G$ on a Hilbert space $\mathcal{H}_\pi.$ For $u,v\in\mathcal{H}_\pi,$ let $\pi_{u,v}$ denote the coefficient function corresponding to $\pi,u$ and $v.$ The {\it Fourier-Stieltjes algebra} of $G,$ denoted $B(G),$ is defined as the collection of all coefficient functions arising from all the unitary representations. In \cite{Eym1}, Eymard introduced the algebra $B(G).$ He showed that it is also the dual of the group $C^*$-algebra $C^*(G).$ With the dual norm $B(G)$ becomes a commutative Banach algebra with the pointwise addition and multiplication.
		
		The closed linear span of all coefficient functions arising only from the left regular representation, $\lambda,$ is called as the {\it Fourier algebra} of the group $G,$ denoted $A(G).$ In \cite{Eym1}, it is proved that $A(G)$ is also a commutative Banach algebra. Further, $A(G)$ is also a regular, semisimple Banach algebra with the Gelfand spectrum homeomorphic to $G.$ It is also a closed ideal of $B(G).$ In the case when $G$ is compact, every irreducible unitary representation  of $G$ is unitarily equivalent to a subrepresentation of the left regular representation. Hence, in this case, $A(G)=B(G).$
		
		For more on the Fourier and the Fourier-Stieltjes algebra, we refer the readers to the fundamental paper of Eymard \cite{Eym1}. Also see \cite{KL2}.
		
	    An {\it $n$-measure} on a locally compact group $G$ is a completely bounded multilinear functional on $C^\ast(G)^n,$ where $C^\ast(G)^n=C^\ast(G)\times C^\ast(G)\times\cdots C^\ast(G)\ (n\mbox{-times}).$ Since $B(G)$ coincides with the dual of $C^\ast(G),$ it follows that the space of all completely bounded multilinear functionals on $C^\ast(G)^n$ coincides with $\otimes_n^{eh}B(G).$ When $G$ is abelian and $n=2,$ this space coincides with the space of all bimeasures on $\widehat{G}$ \cite{GS}. Thus the {\it multidimensional Fourier-Stieltjes algebra}, denoted $B^n(G),$ is the space $\otimes_n^{eh}B(G).$ By \cite{ER2}, $B^n(G)$ coincides with the dual of $\otimes_n^{h}C^\ast(G)$ completely isometrically. With this norm and with pointwise addition and multiplication, $B^n(G)$ becomes a commutative Banach algebra. Further, functions in $B^n(G)$ can be represented as $$f(x_1,x_2,\ldots,x_n)=\langle \pi_1(x_1)\pi_2(x_2)\ldots\pi_n(x_n)\xi,\eta \rangle,$$ where $\pi_i$'s are continuous unitary representations of $G$ on some Hilbert space $\mathcal{H}$ and $\xi,\eta\in\mathcal{H}.$ Also, $$\|f\|_{B^n(G)}=\inf\|\xi\|\|\eta\|,$$ where the infimum is taken over all such representations. In fact, the infimum is actually attained.
		
		Motivated by the definition of the multidimensional Fourier-Stieltjes algebra, the {\it multidimensional Fourier algebra}, denoted $A^n(G),$ is defined as the collection of all functions $f\in L^\infty(G^n)$ such that $\exists$ a normal completely bounded multilinear functional $\Phi$ on $VN(G)^n$ satisfying $$f(x_1,x_2,\ldots,x_n)=\Phi(\lambda(x_1),\lambda(x_2),\ldots\lambda(x_n)).$$ Thus, by \cite{ER2}, $A^n(G)$ coincides with $\otimes^{eh}_n A(G).$ Again, by \cite{ER2} and the fact that the dual of $A(G)$ is $VN(G),$ we have $A^n(G)^\ast$ is completely isometrically isomorphic to $VN^n(G):=\otimes_n^{\sigma h}VN(G).$ If $G$ is compact, then $B^n(G)=A^n(G).$

        For more on the multidimensional Fourier and Fourier-Stieltjes algebra, we refer the readers to the paper of Todorov and Turowska \cite{TT}.

        Throughout the remaining parts of this paper, $G$ from now on will denote a compact group with $dx$ as the normalized Haar measure. Also, we shall denote by $G^n$ the group $G\times G\times\cdots\times G\ (n\mbox{-times}).$
			
\section{Multidimensional Varopolous algebra}\label{MVA}
	In this section, we define the multidimensional version of the Varopolous algebra. For $n=1,$ this algbera, as we know, was defined by Varopolous, assuming that the underlying group is compact and abelian. Here we show that the multidimensional version is a commutative semisimple Banach algebra. 
		
	Let $C(G)$ denote the commutative C*-algebra of continuous complex-valued functions on $G.$ We shall denote by $V^n(G)$ the space $\otimes^h_{n+1}\ C(G).$ If $w\in V^n(G),$ then $w$ can be represented as a norm converging (infinite) sum of elementary tensors and for such $w$ we have $$\|w\|_V=\inf\left\{  \left\| \underset{i=1}{\overset{\infty}{\sum}}\ |\phi^1_i|^2 \right\|^{\frac{1}{2}}_\infty\ldots\left\| \underset{i=1}{\overset{\infty}{\sum}}\ |\phi^{n+1}_i|^2 \right\|^{\frac{1}{2}}_\infty:w=\underset{i=1}{\overset{\infty}{\sum}}\ \phi^1_i\otimes\phi^2_i\otimes\ldots\otimes\phi^{n+1}_i\right\}.$$ For any $\phi\in C(G),$ let $M_\phi$ denote the multiplication operator on $L^2(G).$ Now, for $w=\underset{i=1}{\overset{\infty}{\sum}}\ \phi^1_i\otimes\phi^2_i\otimes\ldots\otimes\phi^{n+1}_i\in V^n(G),$ define the operator $$T_w:\otimes_n^{\sigma h}\ \mathcal{B}(L^2(G))\rightarrow\mathcal{B}(L^2(G))$$ as $$T_w(S_1\otimes S_2\otimes\ldots S_{n})=\underset{i=1}{\overset{\infty}{\sum}}\ M_{\phi^1_i}S_1M_{\phi^2_i}S_2M_{\phi^3_i}\ldots S_nM_{\phi^{n+1}_i}.$$ It is clear that $\|T_w\|=\|w\|_V.$ Also, $T_w$ is weakly continuous in each variable. Further, we have the following.
    \begin{lem}
		If $w\in V^n(G),$ then $T_w$ is completely bounded with $\|T_w\|_{cb}=\|w\|_V$
	\end{lem}
	\begin{proof}
	   Note that, for any $p\in\mathbb{N}$ and $S\in M_p\left( \otimes_n^{\sigma h}\ \mathcal{B}(L^2(G)) \right),$ we have 
		\begin{align*}
			\|T_w^{(p)}(S)\| =& \|M_{\phi_1}S_1^{(p)}M_{\phi_2}S_2^{(p)}M_{\phi_3}\ldots M_{\phi_n}S_{n}^{(p)}M_{\phi_{n+1}}\|\\ \leq& \|\phi_1\|\ldots\|\phi_n\|\|id\otimes S_1\|\ldots\|id\otimes S_n\| \leq \|w\|_V\|S\|,
		\end{align*} 
		i.e., $\|T_w\|_{cb}\leq\|w\|_V.$ For the other way inequality, by Theorem \ref{WASm}, by the fact that the Haagerup tensor product is injective and the inclusion of Haagerup tensor product inside the $\sigma$-Haagerup tensor product, we have $$\|T_\omega\|_{cb}\geq\|T_\omega|_{\otimes_n^h\mathcal{K}(L^2(G))}\|_{cb}= \|\omega\|_V.$$ Hence the proof.
	\end{proof}
	We now show that $V^n(G)$ is a Banach algebra.
	\begin{lem}
		With pointwise addition and pointwise multiplication, the space $V^n(G)$ becomes a commutative Banach algebra.
	\end{lem}
	\begin{proof}
		Note that, it is enough to prove the sub-multiplicativity of the norm as others are clear. In order to prove this, it enough to show the same on the dense subspace $\otimes_{n+1} C(G).$  Observe that, if $$u=\underset{i=1}{\overset{r}{\sum}} \phi_i^1\otimes \phi_i^2\otimes \ldots \otimes \phi_i^{n+1}$$ and $$v=\underset{j=1}{\overset{s}{\sum}} \psi_j^1\otimes \psi_j^2\otimes \ldots \otimes \psi_j^{n+1},$$ then $$u.v=\underset{i=1}{\overset{r}{\sum}}\underset{j=1}{\overset{s}{\sum}}\ \phi_i^1\psi_j^1\otimes \phi_i^2\psi_j^2\otimes \ldots \otimes \phi_i^{n+1}\psi_j^{n+1}$$ and hence $$\|u.v\|_V\leq \left\| \underset{i=1}{\overset{r}{\sum}}\underset{j=1}{\overset{s}{\sum}}\ |\phi_i^1|^2|\psi_j^1|^2 \right\|^{1/2}_\infty\ldots \left\| \underset{i=1}{\overset{r}{\sum}}\underset{j=1}{\overset{s}{\sum}}\ |\phi_i^{n+1}|^2|\psi_j^{n+1}|^2 \right\|^{1/2}_\infty.$$ The last inequality follows from \cite[Proposition 2]{Bl}.
	\end{proof}
	We call the Banach algebra $V^n(G)$ as the multidimensional Varopoulos algebra. By \cite[Pg. 121, Theorem 2.11.2]{Kan1} the Gelfand spectrum of $V^n(G)$ can be identified with $G^{n+1}.$
		
	We now define an action of $G$ on $V^n(G).$ If $x\in G$ and $\omega\in V^n(G),$ then define $x.\omega$ as $$(x.\omega)(x_1,x_2,\ldots,x_{n+1}):=\omega(x_1x,x_2x,\ldots,x_{n+1}x),$$ for $x_1,x_2,\ldots,x_{n+1}\in G.$ It follows as a consequence of the norm on $V^n(G)$ that $(x,\omega)\mapsto x.\omega$ is a continuous action of $G$ on $V^n(G)$ by isometries. Further, as is well known, this action extends to an action of $L^1(G)$ on $V^n(G)$ given by $$f\cdot\omega:=\int_G\ f(x)(x.\omega)\ dx,\ f\in L^1(G)\mbox{ and }\omega\in V^n(G).$$ Note that the above integral is a $V^n(G)$-valued Bochner integral. Further, if $\{f_\alpha\}$ is a bounded approximate identity in $L^1(G)$ and $\omega\in V^n(G),$ then $f_\alpha\cdot\omega\rightarrow\omega.$ As a consequence of Cohen's factorization theorem, we have the following immediate corollary.
	\begin{cor}
		The Banach algebra $V^n(G)$ is an essential $L^1(G)$-module.
	\end{cor}
		%Note that, if $f\in C(G),$ then the mapping $x\mapsto f(x)x\cdot\omega$ is continuous and hence Bochner integrable. Further, $$\|f\cdot\omega\|_{V^n(G)}\leq\|f\|_1\|\omega\|_{V^n(G)}\ \forall\ f\in C(G).$$ Using the fact that $C(G)$ is dense in $L^1(G),$ we can unambiguously define $f\cdot\omega$ for any $f\in L^1(G).$

\section{The maps $P^n$ and $N^n$}\label{PNMap}
	In this section we define the maps $P^n$ and $N^n.$ These maps were originally defined by Varopolous \cite{Var} when the group $G$ is compact and abelian and $n=1.$ Our approach here is motivated by the work of Spronk and Turowska \cite{SprT}.
		
	Let $$V_{inv}^n(G)=\left\{\omega\in V^n(G):\begin{array}{c}\omega(x_1,x_2,\ldots,x_{n+1})=\omega(x_1x,x_2x,\ldots,x_{n+1}x) \\ \forall\ x,x_1,x_2,\ldots,x_{n+1}\in G\end{array}\right\}.$$
	\begin{prop}
		The mapping $P^n:V^n(G)\rightarrow V^n(G)$ given by $$P^n(\omega)(x_1,x_2,\ldots,x_{n+1})=\int_G\omega(x_1x,x_2x,\ldots,x_{n+1}x)\ dx$$ is a contractive projection onto $V_{inv}^n(G)$ and a $V_{inv}^n(G)$-module map.
	\end{prop}
	\begin{proof}
		Since $G$ is compact, the constant function $1$ belongs to $L^1(G)$ and it is clear that $P^n(\omega)=1\cdot\omega\ \forall\ \omega\in V^n(G).$ Hence the map $P^n$ is a contraction. Using the fact that the Haar measure is left invariant, we have that the range of $P^n$ sits inside $V_{inv}^n(G).$ Now, using the definition of $V^n_{inv}(G),$ it is clear that $P^n|_{V^n_{inv}(G}$  is just the identity mapping. This observation also tells us that the map $P^n$ is a $V^n_{inv}(G)$-module map.
	\end{proof}
	Our next theorem is the multidimensional analogue of \cite[Proposition 2.1]{SprT}. This result is motivated by \cite[Theorem 5.7]{TT}.
	\begin{prop}\label{EquiCondUnitMDA}
		Let $u$ be any function on $G^n.$ Then TFAE:
		\begin{enumerate}[(i)]
			\item $u\in b_1(A^n(G))$
			\item $\exists$ an operator $M_u\in b_1\left( \mathcal{B}^\sigma\left( \otimes_n^{\sigma h}\ VN(G), VN(G) \right) \right)$ such that $$M_u(\lambda(x_1)\otimes\ldots\lambda(x_n))=u(x_1,\ldots,x_n)\lambda(x_1\ldots x_n);$$
			\item $\exists$ an operator $\widetilde{M_u}\in b_1\left( \mathcal{B}^\sigma\left( \otimes_n^h\ C^\ast(G), C^\ast(G) \right) \right)$ such that $$\widetilde{M_u}(\lambda(f_1)\otimes\lambda(f_2)\otimes\ldots\otimes\lambda(f_n))=\lambda(g),$$ where $$g(x)=\underset{G^{n-1}}{\int}f_1(x_1)f_2(x_1^{-1}x_2)\ldots f_n(x_{n-1}^{-1}x)u(x_1,x_1^{-1}x_2,\ldots,x_{n-1}^{-1}x)\ dx_1\ldots dx_{n-1}.$$ 
		\end{enumerate}
	\end{prop}
	\begin{proof}
		{\bf (i)$\Rightarrow$(ii)}. Define $\theta:A(G)\rightarrow A^n(G)$ as $$\theta(v)(x_1,x_2,\ldots,x_n):=v(x_1x_2\ldots x_n).$$ By \cite[Proposition 5.2]{TT}, $\theta$ is a complete contraction. Define $m_u:A(G)\rightarrow A^n(G)$ as $m_u(v):=u\theta(v).$ Now, for any $v\in A^n(G)$, we have 
		\begin{eqnarray*}
			\|m_u(v)\|_{A^n(G)} &=& \|u\theta(v)\|_{A^n(G)} \leq \|u\|_{A^n(G)}\|\theta(v)\|_{A^n(G)} \\ &<& \|\theta(v)\|_{A^n(G)} \leq \|v\|_{A^n(G)},
		\end{eqnarray*}
		i.e., $\|m_u\|\leq 1$ or in other symbols $m_u\in b_1(\mathcal{B}(A(G),A^n(G))).$ Let $M_u$ denote the Banach space adjoint of $m_u.$ Then $M_u$ satisfies $(ii).$
	
		{\bf (ii)$\Rightarrow$(iii)}. Let $\widetilde{M_u}=M_u|_{\otimes_n^h C^\ast(G)}.$ Then, by \cite[Pg. 19]{TT}, (iii) follows.
			
		{\bf (iii)$\Rightarrow$(i)}. By \cite{TT}, the spaces $A^n(G)$ and $\left( \otimes_n^h C^\ast(G) \right)^\ast$ are completely isometric to each other. Therefore, let $m_u$ be the Banach space adjoint of $\widetilde{M_u}.$ Then $m_u$ satisfies (i).
	\end{proof}
    Let us now state our next result, which is the multidimensional analogue       of \cite[Theorem 2.2]{SprT}.
	\begin{thm}
		The mapping $N^n:A^n(G)\rightarrow V_{inv}^n(G)$ given by $u\mapsto\widetilde{u},$  where $$\widetilde{u}(x_1,x_2,\ldots,x_{n+1}) :=u(x_1x_2^{-1},x_2x_3^{-1},\ldots,x_nx_{n+1}^{-1}),$$ is an isometric isomorphism.
	\end{thm}
	\begin{proof}
		Let $u\in A^n(G).$ By \cite[Theorem 4.1]{TT} and the fact that $G$ is compact, $u$ can be written as $$u(x_1,x_2,\ldots,x_n)=\left\langle \lambda(x_1)\ldots\lambda(x_n)f,g \right\rangle,$$ for some $f,g\in L^2(G).$ In fact, $f$ and $g$ can be chosen such that $$\|u\|_{A^n(G)}=\|f\|_2\|g\|_2.$$ Using the above form of $u,$ it can be easily shown that $$N^nu(x_1,x_2,\ldots,x_{n+1})=\left\langle \lambda(x_{n+1})^\ast f,\lambda(x_1)^\ast g \right\rangle.$$ For any $\pi\in\widehat{G}$ and $1\leq i\leq d_\pi^2,$ let $$\phi_{\pi,i}(s)=\langle e_i^\pi , \lambda(s)^\ast g \rangle\mbox{ and }\psi_{\pi,i}(s)=\langle \lambda(s)^\ast f, e_i^\pi \rangle,\ s\in G.$$ Note that only countably many of the above are non-zero and hence by Parseval's formula we have $$N^nu(x_1,x_2,\ldots,x_n)=\underset{\pi\in\widehat{G}}{\sum}\underset{i=1}{\overset{d_\pi^2}{\sum}}\phi_{\pi,i}(x_1)\psi_{\pi,i}(x_{n+1}).$$ It follows as in the proof of \cite[Theorem 2.2]{SprT} that both  $$\underset{\pi\in\widehat{G}}{\sum}\underset{i=1}{\overset{d_\pi^2}{\sum}}|\phi_{\pi,i}|^2\mbox{ and } \underset{\pi\in\widehat{G}}{\sum}\underset{i=1}{\overset{d_\pi^2}{\sum}}|\psi_{\pi,i}|^2$$ converge uniformly in $C(G).$ In particular, $N^nu\in V^n(G).$ Further, it can be easily seen that, for any $x\in G,$ $$N^nu(x_1x,x_2x,\ldots,x_{n+1}x) = N^nu(x_1,x_2,\ldots,x_{n+1}),$$ i.e., $N^nu\in V^n_{inv}(G).$ Also, $$\|N^nu\|_V\leq \left\| \underset{\pi\in\widehat{G}}{\sum}\underset{i=1}{\overset{d_\pi^2}{\sum}}|\phi_{\pi,i}|^2 \right\|_\infty^{1/2} \left\|\underset{\pi\in\widehat{G}}{\sum}\underset{i=1}{\overset{d_\pi^2}{\sum}}|\psi_{\pi,i}|^2\right\|_\infty^{1/2} = \|f\|_2\|g\|_2=\|u\|_{A^n(G)}.$$ 
			
		Now, given $\omega\in V^n_{inv}(G),$ define  $$u_\omega(x_1,x_2,\ldots,x_n):=\omega(e,x_1^{-1},x_2^{-1}x_1^{-1}, \ldots,x_n^{-1}x_{n-1}^{-1}\ldots x_1^{-1}).$$ Then, it is a routine check that $N^nu_\omega=\omega,$ thereby showing that the map $N^n$ is surjective. Thus, the only thing left to finish the proof is to show that $u_\omega\in A^n(G)$ and $$\|u_\omega\|_{A^n(G)}\leq\|\omega\|_V.$$ Let $\omega=\underset{i=1}{\overset{\infty}{\sum}}\phi^1_i\otimes\ldots\phi^{n+1}_i.$ Then, for $x_1,x_2,\ldots,x_{n+1}\in G$ and $f\in L^2(G),$ we have,
		\begin{eqnarray*}
			& & T_\omega(\lambda(x_1)\otimes\ldots\lambda(x_n))(f)(x_{n+1}) = \underset{i=1}{\overset{\infty}{\sum}} M_{\phi_i^1}\lambda(x_1)M_{\phi_i^{2}}\ldots M_{\phi_i^{n+1}}f(x_{n+1}) \\ &=& \underset{i=1}{\overset{\infty}{\sum}} \phi^1_i(x_{n+1})\phi^2_i(x_1^{-1}x_{n+1})\ldots\phi_i^{n+1}(x_n^{-1}\ldots x_1^{-1}x_{n+1})f(x_n^{-1}\ldots x_1^{-1}x_{n+1}) \\ &=& \omega(x_{n+1},x_1^{-1}x_{n+1},\ldots x_n^{-1}\ldots x_1^{-1}x_{n+1})\lambda(x_1x_2\ldots x_n)f(x_{n+1}) \\ &=& \omega(e,x_1^{-1},\ldots,x_n^{-1}\ldots x_1^{-1}) \lambda(x_1x_2\ldots x_n)f(x_{n+1}) \\ &=& u_\omega(x_1,x_2,\ldots,x_n)\lambda(x_1x_2\ldots x_n)f(x_{n+1}).
		\end{eqnarray*}
		Since $T_\omega\in\mathcal{B}^\sigma\left( \otimes_n^{\sigma h}\mathcal{B}(L^2(G)),\mathcal{B}(L^2(G) \right),$ we have $T_\omega(\otimes_n^{\sigma h}VN(G))\subseteq VN(G).$ Let $M_{u_\omega}=T_\omega|_{\otimes_n^{\sigma h}VN(G)}.$ Then $M_{u_\omega}$ satisfies (ii) of Proposition \ref{EquiCondUnitMDA} and hence $\|M_{u_\omega}\|\leq\|T_\omega\|=\|\omega\|_V,$ which in turn tells us that $u_\omega\in A^n(G)$ along with the required inequality.
	\end{proof}
	Before we end this section, we define one more map, in terms of $P^n$ and $N^n.$ Define $Q^n:V^n(G)\rightarrow A^n(G)$ by $Q^n=(N^n)^{-1}\circ P^n.$ Therefore, if $\omega\in V^n(G),$ then $$Q^n\omega(x_1,x_2,\ldots,x_n)=\int_G\omega(x,x_1^{-1}x,x_2^{-1}x_1^{-1}x,\ldots,x_n^{-1}x_{n-1}^{-1}\ldots x_1^{-1}x)\ dx.$$ Also, $Q^n,$  being a compotion of two continuous functions is also continuous.
		
\section{Spectral Synthesis in $A^n(G)$ and $V^n(G)$}\label{SSAV}
    In this final and main section, we prove the relation between spectral synthesis in $A^n(G)$ and spectral synthesis in $V^n(G).$ We also prove the relation between strong Ditkin sets in $A^n(G)$ and strong Ditkin sets in $V^n(G).$ 

    \subsection{Submodules of $A^n(G)^\ast$ and $V^n(G)^\ast$}\label{SMAV}
	   We first lay the foundations for the main results of this section. Here we give a correspondence between $A^n(G)$-submodules of $A^n(G)^\ast$ and $V^n(G)$-submodules of $V^n(G)^\ast.$ For the results corresponding to $n=1,$ see \cite{PP1}.
	
	   For an $A^n(G)$-submodule $X$ of $VN^n(G),$ define $$X_{V^n}:=\{S\in V^n(G)^\ast:(\omega\cdot S)\circ (N^n)\in X\ \forall\ \omega\in V^n(G)\}.$$ It is clear that $X_{V^n}$ is a $V^n(G)$-submodule of $V^n(G)^\ast.$ It is also clear that $X_{V^n}$ is w*-closed whenever $X$ is so.
	
	   On the other hand, given a $V^n(G)$-submodule $Y$ of $V^n(G)^\ast,$ let $$Y_{A^n}=\{T\in VN^n(G):(u\cdot T)\circ (N^n)^{-1}\circ P^n\in Y\ \forall\ u\in A^n(G)\}.$$ As expected, $Y_{A^n}$ is a $A^n(G)$-submodule of $VN^n(G)$ and is also w*-closed whenever $Y$ is so.
	
		We now show that the above correspondence is well behaved. For this we need a lemma.
		\begin{lem}\label{Realtion1}
			Let $\omega\in VN^n(G)$ and $T\in VN^n(G).$ Then $$\omega\cdot(T\circ(N^n)^{-1}\circ P^n)\circ N^n=(N^n)^{-1}P^n\omega\cdot T.$$
		\end{lem}
		\begin{proof}
			For any $u\in A^n(G),$ using the fact that $P^n$ is a $V_{inv}^n(G)$-module map, we have,
			\begin{align*}
				\langle u, \omega\cdot(T\circ(N^n)^{-1}\circ P^n)\circ N^n \rangle =& \langle N^nu, \omega\cdot(T\circ(N^n)^{-1}\circ P^n) \rangle \\ =& \langle \omega N^nu, T\circ(N^n)^{-1}\circ P^n \rangle \\ =& \langle P^n(\omega) N^nu, T\circ(N^n)^{-1} \rangle \\ =& \langle N^n((N^n)^{-1}(P^n\omega)u) , T\circ(N^n)^{-1} \rangle \\ =& \langle (N^n)^{-1}(P^n\omega)u , T \rangle \\ =& \langle u , (N^n)^{-1}(P^n\omega)\cdot T \rangle. \qedhere
			\end{align*}
		\end{proof}
		\begin{thm}
			Let $X$ be an $A^n(G)$-submodule of $VN^n(G).$ Then $(X_{V^n})_{A^n}=X.$
		\end{thm}
		\begin{proof}
			Assume that $T\in (X_{V^n})_{A^n}.$ Then, $(u\cdot T)\circ(N^n)^{-1}\circ P^n\in X_{V^n},$ for all $u\in A^n(G),$ which, in turn, implies that $\omega\cdot((u\cdot T)\circ(N^n)^{-1}\circ P^n)\circ N^n\in X,$ for all $\omega\in V^n(G).$ Choosing $u=1$ and $\omega=1\otimes 1\in V^n(G)$ and using Lemma \ref{Realtion1}, we get $T\in X.$
		
			For the converse, let us assume that $T\in X.$ We claim that $(u\cdot T)\circ(N^n)^{-1}\circ P^n\in X_{V^n}.$ But, this is equivalent to showing that $\omega\cdot((u\cdot T)\circ(N^n)^{-1}\circ P^n)\circ N^n\in X$ for all $u\in A^n(G)$ and for all $\omega\in V^n(G).$ But this is a consequence of Lemma \ref{Realtion1} and the fact that $X$ is an $A^n(G)$-submodule.
		\end{proof}
  
		\subsection{Spectral synthesis in commutative Banach algebras}
            Let $\mathcal{A}$ be a regular, semisimple, commutative Banach algebra with the Gelfand structure space $\Delta(\mathcal{A}).$ For a closed ideal $I$ of $\mathcal{A},$ the zero set of $I,$ denoted $Z(I),$ is a closed subset of $\Delta(\mathcal{A})$ defined as $$Z(I)=\{x\in E:\widehat{a}(x)=0\ \forall\ x\in E\}.$$ For a closed subset $E\subset\Delta(\mathcal{A}),$ we define the following ideals in $\mathcal{A}:$
		\begin{eqnarray*}
			j_{\mathcal{A}}(E) &=& \{a\in\mathcal{A}:\widehat{a}\mbox{ has compact support disjoint from E} \}\\
			J_{\mathcal{A}}(E) &=& \overline{j_{\mathcal{A}}(E)} \\
			I_{\mathcal{A}}(E) &=& \{a\in\mathcal{A}: \widehat{a} = 0\mbox{ on }E\}.
		\end{eqnarray*}
		Note that $J_{\mathcal{A}}(E)$ and $I_{\mathcal{A}}(E)$ are closed ideals in $\mathcal{A}$ with the zero set equal to $E$ and $j_{\mathcal{A}}(E)\subseteq I \subseteq I_{\mathcal{A}}(E)$ for any ideal $I$ with zero set $E.$ $E$ is said to be a {\it set of spectral synthesis} (or a \textit{spectral set}) for $\mathcal{A}$ if $I_{\mathcal{A}}(E) = J_{\mathcal{A}}(E).$ 
		
		$E$ is a \textit{Ditkin set} if for every $a\in I_{\mathcal{A}}(E),$ there exists a sequence $\{a_n\}\subset j_{\mathcal{A}}(E)$ such that $u.u_n$ converges in norm to $a.$ If the sequence can be chosen in such a way that it is bounded and is the same for all $a\in I_{\mathcal{A}}(E),$ then we say that $E$ is a \textit{strong Ditkin set.} Note that every Ditkin set is a set of spectral synthesis. 
		
		For more on spectral synthesis see \cite{Kan1,Rei}.

   		Let $X$ be a $\mathcal{A}$-submodule of $\mathcal{A}^\ast.$ For a closed subset $E\subset\Delta(\mathcal{A}),$ we define the following ideals in $\mathcal{A}:$
		\begin{eqnarray*}
			I_{\mathcal{A}}^X(E) &=& \{a\in\mathcal{A}:\langle a,\varphi \rangle = 0\ \forall\ \varphi\in X\cap I_{\mathcal{A}}(E)\} \\
			J_{\mathcal{A}}^X(E) &=& \{a\in\mathcal{A}:\langle a,\varphi \rangle = 0\ \forall\ \varphi\in X\cap J_{\mathcal{A}}(E)\}.
		\end{eqnarray*}
		We say that $E$ is a $X$-spectral set if $\varphi\in X$ with $supp(\varphi)\subseteq E,$ then $\varphi\in I_\mathcal{A}(E)^\perp.$ It is shown in \cite[Proposition 2.4]{PP1} that $E$ is a set of $X$-synthesis if $I_{\mathcal{A}}^X(E)=J_{\mathcal{A}}^X(E).$

	\subsection{Main result}	
		For a closed subset $E$ of $G^n,$ let $E^*$ denote the subset of $G^{n+1}$ given by $$E^\ast:=\{(x_1,x_2,\ldots,x_{n+1})\in G^{n+1}:(x_1x_2^{-1},x_2x_3^{-1},\ldots,x_nx_{n+1}^{-1})\in E\}.$$
		\begin{lem}\label{RelIdeals1}
			Let $E$ be a closed subset of $G^n$ and let $u\in A^n(G).$
			\begin{enumerate}[(i)]
				\item $u\in I_A(E)$ if and only if $N^nu\in I_V(E^\ast).$
				\item $u\in J_A(E)$ if and only if $N^nu\in J_V(E^\ast).$
			\end{enumerate}
		\end{lem}
		\begin{proof}
			Note that (i) is an immediate consequence of the definitions of $N^n$ and $E^\ast.$ To prove (ii), consider the map $\theta^n:G^{n+1} \rightarrow G^n$ given by $\theta^n(x_1,x_2,\ldots,x_{n+1}) = (x_1x_2^{-1},x_2x_3^{-1},\ldots,x_nx_{n+1}^{-1}).$ It can be easily seen that $\theta^n$ is open and continuous. Further, $N^nu=u\circ \theta^n\ \forall\ u\in A^n(G).$ Also, for $u\in A^n(G),$ it can be observed that $supp(N^nu)\cap E^*=(\theta^n)^{-1}(supp(u)\cap E).$ Thus, $u\in j_A(E)$ if and only if $N^nu\in j_V(E^\ast).$ Due to the continuity of the map $N^n,$ we have that $N^nu\in J_V(E^\ast)$ if $u\in J_A(E).$ 
			
			For the other way inclusion, suppose that $N^nu\in J_V(E^*).$ Let $\{\omega_n\}$ be a sequence in $j_V(E^*)$ converging to $N^nu.$ Then, by the same observation made above, we have $supp(N^nQ^n\omega_n)\cap E^\ast=(\theta^n)^{-1}(supp(Q^n\omega)\cap E).$ Also, it can be easily seen that $supp(Q^n\omega_n)\subseteq \theta^n(supp(\omega_n)).$ These two observations, leads to the fact that $Q^n\omega_n\in j_A(E).$ Now, using the continuity of the map $Q^n,$ we have, $$u=Q^nN^nu=\lim Q^n\omega_n\in J_A(E),$$ thereby proving the other way inclusion and hence the lemma.
		\end{proof}
		\begin{lem}
			Let $E\subseteq G^n$ be closed and let $X=J_{A}(E)^\perp=\{T\in VN^n(G):supp(T)\subseteq E\}.$ Then $X_{V^n}=J_V(E^\ast)^\perp=\{S\in V^n(G)^\ast:supp(S)\subseteq E^\ast\}.$
		\end{lem}
		\begin{proof}
			Assume that $S\in X_{V^n}.$ Let $(x_1,x_2,\ldots,x_{n+1})\in supp(S).$ Let $u\in A^n(G)$ such that $u(x_1x_2^{-1},x_2x_3^{-1},\ldots,x_nx_{n+1}^{-1})\neq 0.$ Then $N^nu(x_1,x_2,\ldots,x_{n+1})\neq 0,$ which in turn implies that $(N^nu)\cdot S\neq 0.$ Therefore, there exists $\omega\in V^n(G)$ such that $\langle \omega,(N^nu)\cdot S \rangle\neq 0.$ But $$\langle \omega,(N^nu)\cdot S \rangle = \langle N^nu,\omega\cdot S \rangle = \langle u,(\omega\cdot S)\circ N^n \rangle = \langle 1,u\cdot((\omega\cdot S)\circ N^n) \rangle.$$ Thus, $\langle 1,u\cdot((\omega\cdot S)\circ N^n) \rangle \neq 0,$ i.e., $u\cdot((\omega\cdot S)\circ N^n)\neq 0$ or equivalently, $(x_1x_2^{-1},x_2x_3^{-1},\ldots,x_nx_{n+1}^{-1})\in supp((\omega\cdot S)\circ N^n)\subseteq E.$ 
			
			Conversely, let $S\in V^n(G)^\ast$ such that $supp(S)\subseteq E^\ast.$ Let $\omega\in V^n(G)$ and let $u\in J_A(E).$ Since $u\in J_A(E),$ by Lemma \ref{RelIdeals1}, $N^nu\in J_V(E^\ast)$ and hence $\omega N^nu\in J_V(E^\ast).$ Therefore, $$\langle u,(\omega\cdot S)\circ N^n \rangle = \langle N^nu,\omega\cdot S \rangle = \langle \omega\cdot N^nu,S \rangle=0.$$ Since $u$ and $\omega$ are arbitrary, the proof is complete.
		\end{proof}
		Our next lemma generalises Lemma \ref{RelIdeals1}, but it is worth noting that the proof it makes use of the same.
		\begin{lem}\label{RelIdeals2}
			Let $E$ be a closed subset of $G^n$ and let $u\in A^n(G).$
			\begin{enumerate}[(i)]
				\item $u\in I_A^X(E)$ if and only if $N^nu\in I_V^{Y}(E^\ast).$
				\item $u\in J_A^X(E)$ if and only if $N^nu\in J_V^Y(E^\ast).$
			\end{enumerate}
		\end{lem}
		\begin{proof}
			In view of Lemmas \ref{Realtion1} and \ref{RelIdeals1}, the proof of this follows exactly as in \cite[Lemma 5.3]{PP1}.
		\end{proof}
		
		Here is the first promised result on the relation between spectral synthesis in $A^n(G)$ and spectral synthesis in $V^n(G).$ Given a $A^n(G)$-submodule $X$ of $A^n(G)^\ast,$ in this, we shall denote by $Y$ the corresponding $V^n(G)$-submodule of $V^n(G)^\ast.$
		\begin{thm}\label{SSFAVA}
			A closed subset $E$ of $G^n$ is a set of $X$-synthesis for $A^n(G)$ if and only if $E^\ast$ is a set of $Y$-synthesis for $V^n(G).$
		\end{thm}
		\begin{proof}
			Note that the backward part is an easy consequence of Lemma \ref{RelIdeals2}. We shall now prove the forward part. Here we closely follow Spronk and Turowska \cite[Theorem 3.1]{SprT}. Suppose that $E$ is a set of synthesis for $A^n(G).$ Let $\omega\in I_V^Y(E^\ast).$ Consider a continuous irreducible unitary representation $\pi$ of $G$ with $\mathcal{H}_\pi$ as the representation space. Define two functions $\omega^\pi$ and $\widetilde{\omega}^\pi$ on $G^{n+1}$ with values in $\mathcal{H}_\pi,$ as $$\omega^\pi(x_1,x_2,\ldots,x_{n+1})=\int_G \omega(x_1x,x_2x,\ldots,x_{n+1}x)\pi(x)\ dx$$ and $$\widetilde{\omega}^\pi(x_1,x_2,\ldots,x_{n+1})=\pi(x_1)\omega^\pi(x_1,x_2,\ldots,x_{n+1}).$$ Since $\mathcal{H}_\pi$ is finite dimensional, let us fix an orthonormal basis and write $u_{ij}^\pi$ for the coefficient functions of $\pi.$ Also, let $\omega^\pi_{ij}=u_{ij}^\pi\omega^\pi$ and $\widetilde{\omega}^\pi_{ij}=u_{ij}^\pi\widetilde{\omega}^\pi.$ Note that $\omega^\pi_{ij}\in I_V^Y(E^\ast)$ and hence $\widetilde{\omega}^\pi_{ij}=\underset{k=1}{\overset{d_\pi}{\sum}}(u_{ik}^\pi\otimes 1\otimes\ldots\otimes 1)\omega_{kj}^\pi\in I_V^Y(E^\ast).$ Observe that $$\widetilde{\omega}^\pi(x_1x,x_2x,\ldots,x_{n+1}x)=\widetilde{\omega}^\pi(x_1,x_2,\ldots,x_{n+1})$$ and hence $\widetilde{\omega}_{ij}^\pi\in V^n_{inv}(G).$ Thus, by Lemma \ref{RelIdeals2}, $(N^n)^{-1}(\widetilde{\omega}_{ij}^\pi)\in I_A^X(E)$ and by the assumption on $E,$ we have $(N^n)^{-1}(\widetilde{\omega}_{ij}^\pi)\in J_A^X(E).$ Applying Lemma \ref{RelIdeals2} once again, we have that $\widetilde{\omega}_{ij}^\pi\in J_V^Y(E^\ast.)$ Since $$\omega^\pi(x_1,x_2,\ldots,x_{n+1})=\pi(x^{-1})\widetilde{\omega}^\pi(x_1,x_2,\ldots,x_{n+1}),$$ we have $$\omega_{ij}^\pi=\underset{k=1}{\overset{d_\pi}{\sum}}(\check{u}_{ik}^\pi\otimes 1\otimes\ldots\otimes 1)\widetilde{\omega}_{kj}^\pi\in J_V^Y(E^\ast),$$ where $\check{u}(x)=u(x^{-1}).$ Using the fact that $V^n(G)$ is an essential $L^1(G)$-module, the proof is now completed by using an approximate identity argument (as in \cite[Theorem 3.1]{SprT}) to conclude that $\omega\in J_V^Y(E^\ast.)$ Thus $E^\ast$ is a set of spectral synthesis for $V^n(G).$
		\end{proof}
		Here is the second promised result on the relation between strong Ditkin sets in $A^n(G)$ and strong Ditkin sets in $V^n(G).$ We would like to remark here that we are not proving for $X$-Ditkin sets.
		\begin{thm}\label{DFAVA}
			A closed subset $E$ of $G^n$ is a strong Ditkin set for $A^n(G)$ if and only if $E^\ast$ is a strong Ditkin set for $V^n(G).$
		\end{thm}
		\begin{proof}
			Suppose that $E$ is a strong Ditkin set for $A^n(G).$ Let $\{u_m\}$ be a bounded sequence in $j_{A}(E)$ such that $u_mu\rightarrow u$ for all $u\in I_A(E).$ Let $\omega\in I_V(E^\ast).$ We shall now make use of the notations used in the previous theorem. For each $\pi\in\widehat{G},$ let $\omega^\pi,$ $\omega^\pi_{ij}\in I_V(E^\ast)\cap V_{inv}^n(G),$ as in the previous theorem. Then $(N^n)^{-1}\omega^\pi_{ij}\in I_A(E)$ and hence $u_m(N^n)^{-1}\omega^\pi_{ij}\rightarrow (N^n)^{-1}\omega^\pi_{ij}.$ Applying $N^n,$ $N^nu_m\omega^\pi_{ij}\rightarrow \omega^\pi_{ij}.$ As mentioned in the previous theorem, choose an approximate identity $\{f_\alpha\}$ for $L^1(G).$ Thus, for each $\alpha,$ $$N^nu_m(f_\alpha\cdot\omega)\rightarrow f_\alpha\cdot\omega.$$ Now, 
			\begin{eqnarray*}
				\|N^nu_m\omega - \omega\|_V &\leq& \|N^nu_m\omega - N^nu_m(f_\alpha\cdot\omega)\|_V \\ & & + \|N^nu_m(f_\alpha\cdot\omega) - f_\alpha\cdot\omega\|_V + \|f_\alpha\cdot\omega - \omega\|_V \\ &\leq& \|N^nu_m\|_V\|f_\alpha\cdot\omega - \omega\|_V \\ & & + \|N^nu_m(f_\alpha\cdot\omega) - f_\alpha\cdot\omega\|_V + \|f_\alpha\cdot\omega - \omega\|_V \\ &\leq& C\|f_\alpha\cdot\omega - \omega\|_V + \|N^nu_m(f_\alpha\cdot\omega) - f_\alpha\cdot\omega\|_V,
			\end{eqnarray*}
			where $C$ is a bound for the sequence $\{N^nu_m\}.$ Now, fix an $\alpha$ so that the first term is small and for this $\alpha,$ choose $n$ large enough so that the second term is also small. As the sequence $\{N^nu_m\}$ is bounded, $E^\ast$ is a strong Ditkin set.
			
			We now prove the converse. This part, as we will see, easily follows from Lemma \ref{RelIdeals2}. Suppose that $E^\ast$ is a strong Ditkin set for $V^n(G).$ Let $\{\omega_m\}$ be a bounded sequence in $j_V(E^\ast)$ such that $$\|\omega_m\omega-\omega\|_V\rightarrow 0 \mbox{ as }m\rightarrow\infty\mbox{ and }\forall\ \omega\in I_V(E^\ast).$$ Let $u\in I_A(E).$ Then, by Lemma \ref{RelIdeals2}, $N^nu\in I_V(E^\ast)$ and hence $$\|\omega_m N^nu-N^nu\|\rightarrow\mbox{ as }m\rightarrow\infty.$$ Using the fact that $P^n$ is continuous and also a $V_{inv}^n(G)$-module map, we have that $$P^n(\omega_m N^nu)\rightarrow N^nu\mbox{ as }m\rightarrow\infty.$$ Hence $(N^n)^{-1}(\omega_m N^nu)\rightarrow u$ in $A^n(G).$ Observe that the sequence $\{\omega_m\}$ is contained in $j_V(E^\ast)$ and hence the sequence $\{(N^n)^{-1}P^n\omega_m\}$ belongs to $j_A(E^\ast),$ thereby showing that $((N^n)^{-1}P^n\omega_m)u\rightarrow u$ in $J_A(E).$ Thus $E$ is a strong Ditkin set and hence the proof of the theorem.
		\end{proof}
		\begin{rem}
			As we mentioned earlier, our last theorem on the relation between strong Ditkin sets in $A^n(G)$ and $V^n(G)$ was proved for Ditkin sets and not for the case of $X$-Ditkin sets. Thus, through this remark, we would like to propose the problem of extending the above result to the case of strong $X$-Ditkin sets.
		\end{rem}
	
\section*{Acknowledgement}
	This research for the first author was supported by a grant from CSIR with grant number 09/086(1470)/2020-EMR-I. The second author would like to thank the Science and Engineering Board, India, for the core research grant with file no. CRG/2021/003087/MS.

\end{document}